\documentclass[12pt]{article}
\usepackage{graphics,graphicx}
\usepackage{amssymb}
\usepackage{amsmath}
\usepackage{amsthm}
\DeclareGraphicsExtensions{.png}
\usepackage{latexsym}
\usepackage{tikz}
\usepackage{xcolor}
\newtheorem{theorem}{Theorem}

\newtheorem{example}{Example}
\newtheorem{definition}{Definition}

\date{~}
\title{Towards  optimal control of  systems with backlash}
\author{M.~d.~R. de Pinho\textsuperscript{a},  M.~Margarida A.~Ferreira\\[1mm] \textsuperscript{a} and Georgi Smirnov\textsuperscript{b}\\[2mm]
\textsuperscript{a}Universidade do Porto, Faculdade de Engenharaia, \\[2mm]
DEEC, SYSTEC, Porto, Portugal.  \\[2mm]
\textsuperscript{b}Universidade do Minho, Dep. Matem\'{a}tica, \\[2mm]
Physics Center of Minho and Porto Universities (CF-UM-UP), \\[2mm]
Campus de Gualtar, Braga, Portugal.\\[2mm]
 mrpinho@fe.up.pt , mmf@fe.up.pt, smirnov@math.uminho.pt }

\begin{document}

\maketitle

\begin{abstract}
In this paper we consider time-optimal control problems for systems with backlash. Such systems are described by second order differential equations coupled  with restrictions modeling the inelastic shocks.  A main feature of such systems is the lack of uniqueness of solution to the Cauchy problem. Here,  we introduce approximation systems where the forces during the impact are taken into account. Such approximations are relevant for two reasons. Firstly, we define a set of solutions as limits of  the solutions to the approximation systems. This set   may be  smaller than the set of  of the solutions  usually considered in the literature. Secondly, such approximations are adequate   to derive necessary condition to the time optimal control of interest. To the best of our knowledge, this is the first attempt to derive necessary conditions of optimality for optimal control problems involving systems with backlash.\end{abstract}


\section{Introduction}

Backlash is present in many mechanical systems. 
The study of systems with backlash mainly concerns stabilization problems \cite{Backlash}.
In this paper we consider time-optimal control problems for systems with backlash. Informally, such systems have the  form:
$$
 \ddot{z}=F(z,\dot{z},u)-N_C(z),\;\; z\in C,\;\; u\in U(t),
$$
where $C$ is a set and $N_C(z)$ is a normal force that does not allow the system to cross the boundary of $C$. The shock with the boundary of $C$ is assumed to be absolutely inelastic, i.e., after the collision the velocity $\dot{z}$ loses its component normal to $\partial C$. 
When  $C=\{ z~:~\psi (z)\leq 0\}$ and  $\psi$ is a smooth function,  $N_C(z)$ has the form $N\nabla\psi(z)$ for some  $N\geq 0$. The study of these systems is hard and it has a long story.
To our knowledge  it was  Moreau \cite{Mor} who first studied systems with completely inelastic impacts. Later,  for  second order differential equations Schatzman introduced the definition of solution  involving   measures with respect to time \cite{Sch}  and Monteiro-Marques \cite{MM}    proved the  first existence result for the Cauchy problem 
\begin{equation}\label{CP}
\begin{array}{lcl}
d\dot{z} & = & F(z,\dot{z},u)dt-\nabla\psi(z)d\nu,\\
z(0) & = & z_0,\\
\dot z(0) & = & \dot z_0,
\end{array}
\end{equation}
where $\nu\in BV([0,T],R^n)$, $d\nu\geq 0$,  $\psi(z)d\nu=0$ and,  if $\psi(z(t))=0$,  then 
$\frac{d}{dt}\psi(z(t+0))=0$. 
The solution is  a pair $(z,\dot{z})\in AC([0,T],R^n)\times BV([0,T],R^n)$ satisfying the above conditions. We  call  it a  {\em Moreau-Schatzman-Monteiro Marques solution}  (MSMM solution). 
For general theorem on the existence of solution see, e.g., \cite{Bernicot}.  It is of foremost importance to notice that the solution  may not be unique \cite{Ballard}. 

The model with  inelastic shocks is  adequate for the description of systems with backlash when the motion of the bodies is rather slow. Moreover,   for the study of optimal control problems, such model is simpler than the case of elastic shocks.

Noteworthly,   the equation 
$$
d\dot{z} = F(z,\dot{z},u)dt-\nabla\psi(z)d\nu
$$ 
is of interest to describe dynamics of systems with hysteresis \cite{KP} (backlash is an example of hysteresis).


In this paper, we focus on time-optimal control problems for such systems. Our aim is to  derive necessary conditions  of optimality. To do so   we introduce approximating systems where forces during the impact are taken into account.  Considering that the forces between interacting bodies go to infinity,  we define solutions of \eqref{CP} as a set of limits of solutions to the approximating systems. This set of solutions, which we call \emph{backlash solutions}, is  included in the set of MSMM solutions. 

The approximating techniques are in vein of our work on optimal control for systems involving sweeping processes \cite{nosso_2022}. Sweeping systems are first order systems. Optimal control problems  involving those systems are nowadays well studied; see, e.g., \cite{MoCa17,colombo, zeidan2020} and references therein.

The study of  the second order system \eqref{CP} is much more involved than that of sweeping systems. Remarkably, and as we show here, the approximating techniques similar to those in  \cite{nosso_2022}  allow us to derive necessary conditions of optimality for backlash solutions.  Central to our analysis is the fact that   necessary condition of optimality  for  time-optimal control problems involving  the  approximating systems are well known. Passing to the limit (i.e., considering that the acting forces go to infinity), we  get the desired necessary conditions of optimality for backlash solutions. 

%
%
%

To keep the focus on the main difficulties regarding the derivation of necessary conditions, we consider a  time-optimal control problem. For simplicity of presentation, we assume the data smooth. This work may be easily generalized to cover optimal control problems with different costs, additional end point constraints and some nonsmooth data.   

Since we  focus on global solutions to the time-optimal control problem of interest, it is  enough to show that there exists an optimal  solution satisfying the necessary conditions. It is in this vein that we formulate our results here. One could use penalization techniques to show that any other optimal  solution also satisfies the necessary conditions. However, this is not our aim; we concentrate on a solution and not on all of them. 


\vspace{5mm}

Throughout this  paper we denote  the set of real numbers by $R$ and
 the usual  $n$-dimensional  space    of  vectors
$x=(x^1,\ldots ,x^n)$,  where $x^i\in R,\;  i=\overline{1,n}$, by $R^n$. 
 The
inner product of two vectors $x$ and $y$ in $R^n$ is defined by 
$$
\langle x,y\rangle =x^1y^1+\ldots +x^ny^n,
$$
and the norm  of  a  vector  $x\in  R^n$  is    $|x|=\langle x,x\rangle^{1/2}$.  
If  $g: R^p\to  R^q$, $\nabla g$ represents the derivative. The partial derivative in order of a variable $\xi$ is denoted by $\partial_{\xi}$. If $a\in R$, we set $a_+=\max\{a,0\}$.
By $\chi_E$ we denote the characteristic function of a set $E\subset [0,1]$.

 The space  $L_{\infty}([a,b]; R^p)$ (or simply $L_{\infty}$ when the domains are clearly understood) is the Lebesgue space of essentially 
bounded functions  $h:[a,b]\to R^p$.  We say that $h\in BV([a,b]; R^p)$ if $h$ is a function of bounded variation. The space of continuous functions is denoted by $C([a,b]; R^p)$ and the space of absolutely continuous functions  by  $AC([a,b]; R^p)$. 

In this work, we use sequences and subsequences of functions indexed by $\gamma$. To simplify the notation, we write $z_\gamma$ meaning $z_{\gamma_k}$ or, when considering subsequences, 
$z_{\gamma_{{k}_m}}$, i.e., we omit the indexes  unless the appearance of the index is absolutely necessary as  in the case of  construction of diagonal subsequences.

\section{Definition of solution}

Consider system \eqref{CP}. 
In many situations, instead of $z$ we can introduce new coordinates $(x,y)$, where $y=\psi(z)$, and rewrite the system in the following {\em canonical} form:
\begin{eqnarray}
&& \ddot{x}=f(x,y,\dot{x},\dot{y},u),\label{01}\\
&& d\dot{y}=g(x,y,\dot{x},\dot{y},u)dt-d\nu,\label{02}\\
&& y\leq 0,\;\; yd\nu=0,\,\; d\nu\geq  0, ~  u \in U\label{03}\\
&& y(t)=0 \Longrightarrow ~ \dot y(t+0)=0\label{04}
\end{eqnarray}
Here, we focus on  systems having this form.
We start by defining  a solution to system (\ref{01}) - (\ref{04}) suitable to derive necessary conditions of optimality. 

We assume that $U(t)\subset R^m$ is a convex compact set, the map $t\rightarrow U(t)$ is bounded and measurable, and that the right-hand side of the control system has the following structure:
$$
f(x,y,v,w,u)=f_1(x,y,v)+f_2(x,y,v)w+f_3(x,y,v)u
$$
 and 
 $$
 g(x,y,v,w,u)=g_1(x,y,v)+g_2(x,y,v)w+g_3(x,y,v)u.
 $$
Furthermore, we assume that $f$ and $g$ are sufficiently smooth functions and 
\begin{equation}\label{fg:bound}
|f|\leq M(1+|x|+|y|+|v|+|w|)\text{ and } |g|\leq M(1+|x|+|y|+|v|+|w|).
\end{equation}

\vspace{5mm}

For  $a\in R$, we denote  $a_+=\max\{a,0\}=\frac{1}{2}(a+|a|)$. Let $\gamma>0$ and $u_{\gamma}\in U$ be a measurable function. 

Consider now  the Cauchy problem
\begin{eqnarray}
&& \dot{x}_{\gamma}=v_{\gamma},\label{1}\\
&& \dot{y}_{\gamma}=w_{\gamma},\label{2}\\
&& \dot{v}_{\gamma}=f(x_{\gamma},y_{\gamma},v_{\gamma},w_{\gamma},u_{\gamma}),\label{3}\\
&& \dot{w}_{\gamma}=g(x_{\gamma},y_{\gamma},v_{\gamma},w_{\gamma},u_{\gamma})-\gamma (y_{\gamma})_+(w_{\gamma})_+,\label{4}\\
&& x_{\gamma}(0)=x_0,\;y_{\gamma}(0)=y_0,\;v_{\gamma}(0)=v_0,\;w_{\gamma}(0)=w_0.\label{5}
\end{eqnarray}

The force $\gamma (y_{\gamma})_+(w_{\gamma})_+$ appearing in  \eqref{4} is basically Hooke's law with the additional term $(w_{\gamma})_+$ introduced to model the inelastic shock: the reaction force is active only when the systems moves outward  the set $C$.
\vspace{5mm}

\begin{theorem}
Let $\gamma_k\rightarrow\infty$ be a sequence of scalars and let $u_{{\gamma_k}_l}$ be a sequence  of admissible controls converging  in weak* topology of  $L_{\infty} ([0,T],R^n)$ to some function $\hat u(t)\in U(t)$.
For each $k_l$,   let $(x_{{\gamma_k}_l},
y_{{\gamma_k}_l},v_{{\gamma_k}_l},w_{{\gamma_k}_l})$   be a solution to the system (\ref{1})-(\ref{5}). 

Then the sequence  $(x_{{\gamma_k}_l},
y_{{\gamma_k}_l},v_{{\gamma_k}_l},w_{{\gamma_k}_l})$ has a  subsequence converging to a limiting function $(\hat x,\hat y,\hat v,\hat w)\in AC([0,T],R^n)\times AC([0,T],R)\times AC([0,T],R^n)\times BV([0,T],R)$ where the convergence in $AC$ is in the sense of the uniform norm and the convergence in $BV$ is in the sense of the weak* topology. Moreover, $(\hat x,\hat y,\hat v,\hat w,\hat u)$ solves the problem

\begin{eqnarray}
&& \dot{x}=v,\label{1a}\\
&& \dot{y}=w,\label{2a}\\
&& \dot{v}=f(x,y,v,w,u),\label{3a}\\
&& dw=g(x,y,v,w,u)dt-d\nu, ~~d\nu\geq 0,~~ yd\nu=0, \label{4a}\\
&& y\leq 0,~~y(t)=0\Longrightarrow w(t+0)=0,\label{4aa}\\
&& x(0)=x_0,\;y(0)=y_0,\;v(0)=v_0,\;w(0)=w_0.\label{5a}
\end{eqnarray}
\end{theorem}

\begin{proof}

We fix a sequence $\gamma_{k_l}\rightarrow\infty$.  
 From (\ref{1}) - (\ref{4}) and \eqref{fg:bound} we have
$$
\frac{1}{2}\frac{d}{dt}(|x_{{\gamma_k}_l}|^2+|y _{{\gamma_k}_l}|^2+|v_{{\gamma_k}_l} |^2+|w _{{\gamma_k}_l}|^2)
$$
$$
\leq 
M_1+M_2(|x_{{\gamma_k}_l} |^2+|y _{{\gamma_k}_l}|^2+|v_{{\gamma_k}_l} |^2+|w_{{\gamma_k}_l} |^2)- {\gamma_k}_l(y_{{\gamma_k}_l} )_+(w_{{\gamma_k}_l} )_+^2
$$
$$
\leq M_1+M_2(|x_{{\gamma_k}_l} |^2+|y _{{\gamma_k}_l}|^2+|v_{{\gamma_k}_l} |^2+|w_{{\gamma_k}_l} |^2).
$$
Applying the Gronwall inequality we see that $|x_{{\gamma_k}_l }|^2+|y_{{\gamma_k}_l} |^2+|v_{{\gamma_k}_l} |^2+|w_{{\gamma_k}} |^2$ is a bounded function. This implies that $(x_{{\gamma_k}_l} ,y_{{\gamma_k}_l} )$ contains a uniformly converging subsequence. Hence (see (\ref{3})) 
$v_{{\gamma_k}_l}$ also contains a uniformly converging subsequence. From (\ref{4}) we get
$$
w_{{\gamma_k}_l}=w_0+\int_0^Tg_{{\gamma_k}_l}dt-\int_0^T {{\gamma_k}_l}(y_{{\gamma_k}_l})_+(w_{{\gamma_k}_l})_+dt
$$
This implies that the sequence of functions
$$\int _0^t{{\gamma_k}_l}(y_{{\gamma_k}_l})_+(w_{{\gamma_k}_l})_+ds$$
 is bounded in $BV([0,T],R)$. Therefore,  the sequence of measures  
 ${{\gamma_k}_l}(y_{{\gamma_k}_l})_+(w_{{\gamma_k}_l})_+dt$
 contains a weak* converging subsequence to some measure $d\nu$, $\nu \in BV([0,T],R)$ (we do not relabel). From (\ref{4}), we deduce that there exists a constant $K>0$ such that
$$
\displaystyle \int_0^T |\dot w_{\gamma_{k_l}} |dt \leq K.
$$
It then follows that, without loss of generality, $w_{\gamma_{k_l}}$ converges pointwisely to some function $\hat w\in BV$.

Going back to (\ref{1}),  (\ref{2}) and (\ref{3}) we get
$$
\begin{array}{rcl}
x_{\gamma_{k_l}} & = & x_0+\int_0^t v_{\gamma_{k_l}} (s)ds,\\[2mm]
y_{\gamma_{k_l}} & = & y_0+\int_0^t w_{\gamma_{k_l}}(s)ds, \\[2mm]
v_{\gamma_{k_l}} & = & v_0+\int_0^t \big( (f_1)_{\gamma_{k_l}} +(f_2)_{\gamma_{k_l}}w_{k_l} (s)\big)ds\\[2mm]
& & \qquad+ \int_0^t \big( (f_3)_{\gamma_{k_l}} -(\hat f_3)\big) u_{\gamma_{k_l}} (s)ds+\int_0^t (\hat f_3) u_{\gamma_{k_l}} (s)ds,
\end{array}
$$
where $(f_i)_{\gamma_{k_l}}=f_i(x_{\gamma_{k_l}},y_{\gamma_{k_l}},v_{\gamma_{k_l}})$, $i=1,2,3$, and $\hat f_3=f_3(\hat x, \hat y, \hat v)$.  Passing to the limit, we get 
$(\hat x, \hat y, \hat v) \in  AC([0,T],R^n)\times AC([0,T],R)\times AC([0,T],R^n)$  satisfying equations 
(\ref{1a}) -- (\ref{3a}).

We now turn to (\ref{4a}). Take any $w^*\in C([0,T],R)$. Turning to (\ref{4}) we have
$$
\int_0^T w^* \dot w_{\gamma_{k_l}}dt= \int_0^T w^*g_{\gamma_{k_l}} dt-  \int_0^T \gamma_{k_l} (y_{\gamma_{k_l}})_+(w_{\gamma_{k_l}})_+w^* dt,
$$
where $g_{\gamma_{k_l}}=(g_1)_{\gamma_{k_l}}+ (g_2)_{\gamma_{k_l}}w_{\gamma_{k_l}}+(g_3)_{\gamma_{k_l}}u_{\gamma_{k_l}}$, $(g_i)_{\gamma_{k_l}}=g_i(x_{\gamma_{k_l}},y_{\gamma_{k_l}},v_{\gamma_{k_l}})$.
Without loss of generality we deduce that the measures $\dot w_{\gamma_{k_l}}dt$ and  $g_{\gamma_{k_l}} dt$  converge to $d\hat w$ and $d\hat G$. Hence, we have
$$
d\hat w= d\hat G- d\nu.
$$
We now turn to $\hat G$. Take any $z^*\in L^1([0,T];R)$.  We have 
$$
\begin{array}{c}
\int_0^T z^* g_{\gamma_{k_l}}dt= \int_0^T z^*\big((g_1)_{\gamma_{k_l}}+ (g_2)_{\gamma_{k_l}}w_{\gamma_{k_l}}\big)dt \\[2mm]
\qquad \quad
 + \int_0^T z^* \big((g_3)_{\gamma_{k_l}}-\hat g_3\big) u_{\gamma_{k_l}}dt
 + \int_0^T z^* \hat g_3u_{\gamma_{k_l}}dt,
\end{array}
$$
where $(g_i)_{\gamma_{k_l}}=g_i(x_{\gamma_{k_l}},y_{\gamma_{k_l}},v_{\gamma_{k_l}})$, $i=1,2,3$, and $\hat g_3=g_3(\hat x, \hat y, \hat v)$.
Passing to the limit, we conclude that
$$
d\hat G= \hat g_1dt+\hat g_2 \hat w dt +\hat g_3 \hat u dt, 
$$
where $\hat g_i=g_i(\hat x, \hat y, \hat v)$, $i=1,2,3$. 
It is a simple matter to see that we get $ d\nu \geq 0$ and $\hat y d\nu=0$. Therefore we get (\ref{4a}). It remains to see that (\ref{4aa}) holds.


Let us show that $\hat y(t)\leq 0$. We  fix $\gamma={\gamma_k}_l$ and omit the respective index.  Let $t_0$ and $\tau$ be such that $y(t_0)=0$ and 
$y(t_0+\tau)=\max\{ y(t)~:~t\in [0,T]\}$,  when  $y(t)>0$ for  $t\in ]t_0,t_0+\tau ]$ and $y(t)< y(t_0+\tau)$ for $t\in [t_0,t_0+\tau[$. We now consider subintervals  $]t_j, t_{j}+\tau_j[\subset [t_0,t_0+\tau]$, $j=0,1,\ldots$ where  $w(t)>0$, $t\in ]t_j,t_{j}+\tau_j[$, $w(t_j)=0$, $w(t_{j}+\tau_j)=0$, $j=1,2,\ldots$. Clearly,  the set of such subintervals is countable and there exists a $t_*\leq t_0+\tau $ such that $t_j\to t_*$ and $\displaystyle\lim_{j\to +\infty }y(t_j)=y(t_*)$. Thus, we deduce  $t_*=t_0+\tau$ (indeed, if $t_*<t_0+\tau$, then there would be another subinterval where $y$ would increase and, by construction of the subintervals $[t_j,t_j+\tau_j]$, this is not possible). 

We now integrate (\ref{4}). Observing that $\gamma (y)_+(w)_+=\frac{\gamma }{2}\frac{d}{dt}y^2$ whenever $y>0$ and $w>0$, we get 
\begin{eqnarray*}
&& \frac{\gamma}{2}(y^2(t_0+\tau_0)-y^2(t_0))=w(t_0)+\int_{t_0}^{t_0+\tau_0}gdt,\\
&& \frac{\gamma}{2}(y^2(t_1+\tau_1)-y^2(t_1))=\int_{t_1}^{t_1+\tau_1}gdt,\\
&& \cdots\\
&& \frac{\gamma}{2}(y^2(t_{j}+\tau_j)-y^2(t_j))=\int_{t_j}^{t_{j}+\tau_j}gdt,\\
&& \cdots\\
\end{eqnarray*}
Since $y(t_{j}+\tau_j)\geq y(t_{j+1})$,  we get
\begin{eqnarray*}
&& \frac{\gamma}{2}y^2(t_1)\leq w(t_0)+\int_{t_0}^{t_0+\tau_0}gdt,\\
&& \frac{\gamma}{2}(y^2(t_2)-y^2(t_1))\leq\int_{t_1}^{t_1+\tau_1}gdt,\\
&& \cdots\\
&& \frac{\gamma}{2}(y^2(t_{j+1})-y^2(t_j))\leq\int_{t_k}^{t_{j}+\tau_j}gdt,\\
&& \cdots\\
\end{eqnarray*}
Adding these inequalities we obtain
$$
\frac{\gamma}{2}y^2(t_0+\tau)\leq w(t_0)+\int_{t_0}^{t_0+\tau}|g|dt\leq w(t_0)+M_3.
$$
Hence
$$
y(t_0+\tau)\leq \frac{M_4}{\sqrt{\gamma}}\rightarrow 0,\;\; \gamma\rightarrow\infty,
$$
and so we have $\hat y(t)\leq 0$.

Now we prove the inelastic shock condition. Notice that if $y(t_0)=0$, $w(t_0)\geq 0$, then, for small $\theta>0$, we have
$$
0\leq w(t)\leq\bar{M}-\frac{\gamma}{2}y^2(t),\;\; t\in [t_0,t_0+\theta].
$$
Since the solution to the Cauchy problem
$$
\dot{y}=\bar{M}-\frac{\gamma}{2}y^2,\;\; y(t_0)=0,
$$
is given by
$$
y(t)=\sqrt{\frac{2\bar{M}}{\gamma}}\frac{e^{\sqrt{2\bar{M}\gamma}(t-t_0)}-1}{e^{\sqrt{2\bar{M}\gamma}(t-t_0)}+1},
$$
we get
$$
w(t)\leq \bar{M}\left( 1-\left(\frac{e^{\sqrt{2\bar{M}\gamma}(t-t_0)}-1}{e^{\sqrt{2\bar{M}\gamma}(t-t_0)}+1}\right)^2\right)\rightarrow 0,\; \gamma\rightarrow\infty.
$$
Thus we have the inelastic shock condition $w(t_0+0)=0$ satisfied.
\end{proof}

Now we are in a position to define the solution to Cauchy problem (\ref{01}) - (\ref{04}).

\begin{definition}\label{backlash:sol}
We say that   $(x,y,v,w)\in AC([0,T],R^n)\times AC([0,T],R)\times AC([0,T],R^n)\times BV([0,T],R)$ 
 is a {\em backlash  solution} corresponding to a control $u$, if $(x,y,v,w,u)$ 
 is a limit of some sequence of admissible control processes $$(x_{{\gamma_k}_l},y_{{\gamma_k}_l},v_{{\gamma_k}_l},w_{{\gamma_k}_l},u_{{\gamma_k}_l})$$ such that 
 $(x_{{\gamma_k}_l},y_{{\gamma_k}_l},v_{{\gamma_k}_l})$ converges uniformly to $(x,y,v)$, $w_{{\gamma_k}_l}$ converges  to $w$ in the weak* topology of $BV([0,T],R)$  and 
$u_{{\gamma_k}_l}$ converges to $u$ in the weak* topology of $L_{\infty}([0,T],R^m)$, as $l\rightarrow\infty$. 
\end{definition}
In general the solution is not unique.

\subsection*{Comparison with MSMM solution}

In general, the set of MSMM solutions is strictly larger than the set of backlash solutions. Indeed, let us consider a point body moving along a line:
$$
\ddot{y}=u-N,\;\; y\leq 0,\;\; y(0)=\dot{y}(0)=0.
$$
There exists $u(\cdot)\in C^{\infty}$ such that the MSMM solution is not unique (see \cite{Ballard}).  Put $U(t)=\{ u(t)\}$ and choose the sequence $\gamma_k\rightarrow\infty$ such that the sequence of approximating solutions converges. So, we have a unique backlash solution. 

\section{Approximating problem}

To simplify the forthcoming analysis we  now proceed assuming that  $U(t)=U$, where $U$ is a compact and convex set.

Let $(\hat{x},\hat{y},\hat{v},\hat{w},\hat{u})$ be a global optimal solution of the problem
\begin{eqnarray}
&& T\rightarrow\min,\label{0b}\\
&& \dot{x}=v,\label{1b}\\
&& \dot{y}=w,\label{2b}\\
&& \dot{v}=f(x,y,v,w,u),\label{3b}\\
&& dw=g(x,y,v,w,u)dt-d\nu,\;\;  ~~d\nu\geq 0,~~ yd\nu=0, \label{4b}\\
&& y\leq 0,~~y(t)=0\Longrightarrow w(t+0)=0,  \label{4bb}\\
&& x(0)=x_0,\;y(0)=y_0<0,\;v(0)=v_0,\;w(0)=w_0,\label{5b}\\
&& x(T)=x_1,\;y(T)=y_1<0,\;v(T)=v_1,\;w(T)=w_1.\label{6b}
\end{eqnarray} 
Let $\hat{T}$ be the optimal time. 
According to the Definition \ref{backlash:sol} there exists a sequence ${\gamma_{k_l}}\rightarrow\infty$ 
  such that 
\begin{equation}\label{limit}
(\hat{x},\hat{y},\hat{v},\hat{w}, \hat u)=\lim_{l\rightarrow\infty}(x_{\gamma_{k_l}},y_{\gamma_{k_l}},v_{\gamma_{k_l}},w_{\gamma_{k_l}},u_{\gamma_{k_l}}),
\end{equation}
where $(x_{\gamma_{k_l}},y_{\gamma_{k_l}},v_{\gamma_{k_l}},w_{\gamma_{k_l}},u_{\gamma_{k_l}})$ is a control process of the system
\begin{eqnarray}
&& \dot{x}_{\gamma_{k_l}}=v_{\gamma_{k_l}},\label{1c}\\
&& \dot{y}_{\gamma_{k_l}}=w_{\gamma_{k_l}},\label{2c}\\
&& \dot{v}_{\gamma_{k_l}}=f(x_{\gamma_{k_l}},y_{\gamma_{k_l}},v_{\gamma_{k_l}},w_{\gamma_{k_l}},{u}_{\gamma_{k_l}}),\label{3c}\\
&& \dot{w}_{\gamma_{k_l}}=g(x_{\gamma_{k_l}},y_{\gamma_{k_l}},v_{\gamma_{k_l}},w_{\gamma_{k_l}},{u}_{\gamma_{k_l}})-{\gamma_{k_l}} (y_{\gamma_{k_l}})_+(w_{\gamma_{k_l}})_+,\label{4c}\\
&& x_{\gamma_{k_l}}(0)=x_0,\;y_{\gamma_{k_l}}(0)=y_0,\;v_{\gamma_{k_l}}(0)=v_0,\;w_{\gamma_{k_l}}(0)=w_0.\label{5c}
\end{eqnarray}
There exists a  sequence $T_l$ such that $\displaystyle \lim_{l\to \infty} T_l=\hat T$ and 
$$
\lim_{l\to \infty}(x_{\gamma_{k_l}}(T_l),y_{\gamma_{k_l}}(T_l),v_{\gamma_{k_l}}(T_l),w_{\gamma_{k_l}}(T_l))=(x_1,y_1,v_1,w_1).
$$
(Here we use the condition that  $y_1<0$. This implies that  $\hat w$ is continuous in a neighborhood of $\hat T$.)

Now, we consider  the following time-optimal control problem:
\begin{eqnarray}
&& T\rightarrow\min,\label{0d}\\
&& \dot{x}=v,\label{1d}\\
&& \dot{y}=w,\label{2d}\\
&& \dot{v}=f(x,y,v,w,u),\label{3d}\\
&& \dot{w}=g(x,y,v,w,u)-{\gamma_{k_l}}y_+w_+,\label{4d}\\
&& x(0)=x_0,\;y(0)=y_0,\;v(0)=v_0,\;w(0)=w_0,\label{5d}\\
&& x(T)=x_{\gamma_{k_l}}(T_l),\;y(T)=y_{\gamma_{k_l}}(T_l),\;v(T)=v_{\gamma_{k_l}}(T_l),\;w(T)=w_{\gamma_{k_l}}(T_l).\label{6d}
\end{eqnarray} 

This problem has an admissible solution $(x_{\gamma_{k_l}},y_{\gamma_{k_l}},v_{\gamma_{k_l}},w_{\gamma_{k_l}},u_{\gamma_{k_l}})$.
Denote by 
$$(\hat{x}_{\gamma_{k_l}},\hat{y}_{\gamma_{k_l}},\hat{v}_{\gamma_{k_l}},\hat{w}_{\gamma_{k_l}},\hat{u}_{\gamma_{k_l}})$$ an optimal solution to this problem. The corresponding optimal time is $\hat{T}_l\leq T_l$, with $\displaystyle \liminf _{l\to \infty} \hat{T}_l=\hat T$ and, without loss of generality, we can consider that $\hat{T}_{l}$ is a  monotone increasing sequence.

\vspace{5mm}

Applying necessary conditions of optimality to this problem we see that there exist $(q,s,p,r)\in AC([0,\hat{T}],R^{n+1+n+1})$ (we omit the index $\gamma_{k_{l}}$) and a constant $H_l$ satisfying the following conditions (see \cite[Theorem 8.7.1]{Vinter}):

\begin{eqnarray}
&& \dot{q}=-(\nabla_x\hat{f}_{\gamma_{k_l}})^*p-\nabla_x\hat{g}_{\gamma_{k_l}}r,\label{pm1}\\
&&\dot{s}=-\langle \partial_y\hat{f}_{\gamma_{k_l}},p\rangle-\partial_y\hat{g}_{\gamma_{k_l}}r+{\gamma_{k_l}}h_y(\hat{w}_{\gamma_{k_l}})_+r,\label{pm2}\\
&& \dot{p}=-q-(\nabla_v\hat{f}_{\gamma_{k_l}})^*p-\nabla_v\hat{g}_{\gamma_{k_l}}r,\label{pm3}\\
&& \dot{r}=-s-\langle \partial_w\hat{f}_{\gamma_{k_l}},p\rangle-\partial_w\hat{g}_vr+{\gamma_{k_l}}h_w(\hat{y}_{\gamma_{k_l}})_+r,\label{pm4}\\
&& \max_{u\in U}(\langle q, \hat{v}_{\gamma_{k_l}}\rangle+s\hat{w}_{\gamma_{k_l}}+\langle p, 
f(\hat{x}_{\gamma_{k_l}},\hat{y}_{\gamma_{k_l}},\hat{v}_{\gamma_{k_l}},\hat{w}_{\gamma_{k_l}},u)\rangle \nonumber\\
&& \nonumber  \qquad   \qquad +
r(g(\hat{x}_{\gamma_{k_l}},\hat{y}_{\gamma_{k_l}},\hat{v}_{\gamma_{k_l}},\hat{w}_{\gamma_{k_l}},u)-{\gamma_{k_l}}(\hat{y}_{\gamma_{k_l}})_+(\hat{w}_{\gamma_{k_l}})_+))\nonumber\\
&& =\langle q, \hat{v}_{\gamma_{k_l}}\rangle+s\hat{w}_{\gamma_{k_l}}+
\langle p, \hat{f}_{\gamma_{k_l}} \rangle +
r(\hat{g}_{\gamma_{k_l}}-{\gamma_{k_l}}(\hat{y}_{\gamma_{k_l}})_+(\hat{w}_{\gamma_{k_l}})_+)\equiv  H_l \geq 0, \label{pm5}\\[2mm]
&& |q(\hat{T}_l)|^2+|s(\hat{T}_l)|^2+|p(\hat{T}_l)|^2+|r(\hat{T}_l)|^2=1,\label{pm6}
\end{eqnarray}
where $\hat{f}_{\gamma_{k_l}}=f(\hat{x}_{\gamma_{k_l}},\hat{y}_{\gamma_{k_l}},\hat{v}_{\gamma_{k_l}},\hat{w}_{\gamma_{k_l}},\hat{u}_{\gamma_{k_l}})$, 
$\hat{g}_{\gamma_{k_l}}=g(\hat{x}_{\gamma_{k_l}},\hat{y}_{\gamma_{k_l}},\hat{v}_{\gamma_{k_l}},\hat{w}_{\gamma_{k_l}},\hat{u}_{\gamma_{k_l}})$.
The functions
$$
h_y(t)=\left\{
\begin{array}{cl}
0, & \hat{y}_{\gamma_{k_l}}(t)<0,\\
h\in [0,1], & \hat{y}_{\gamma_{k_l}}(t)=0,\\
1, & \hat{y}_{\gamma_{k_l}}(t)>0,
\end{array}
\right.
$$
and
$$
h_w(t)=\left\{
\begin{array}{cl}
0, & \hat{w}_{\gamma_{k_l}}(t)<0,\\
h\in [0,1], & \hat{w}_{\gamma_{k_l}}(t)=0,\\
1, & \hat{w}_{\gamma_{k_l}}(t)>0,
\end{array}
\right.
$$
are measurable.

\section{Main assumption}


Recall that the set of  density points of any measurable set $A\subset [0,\hat{T}]$ is
$$
D_A=\{t\in A~:~\exists\; \lim_{\delta\downarrow 0}\frac{{\rm meas\:}([t-\delta,t+\delta]\cap A)}{2\delta}=1\}.
$$
It is well known that  ${\rm meas\:}(A\setminus D_A)=0$.

\vspace{6mm}


Let $E\subset [0,\hat{T}]$ be the set such that if $t\in E$,
 then  $t$ is a Lebesgue point of 
$\dot{\hat{v}}_{\gamma_{k_l}}$, $\dot{\hat{w}}_{\gamma_{k_l}}$,  $\hat{u}_{\gamma_{k_l}}$, $\dot{q}$, $\dot{s}$, $\dot{p}$, and $\dot{r}$  simultaneously   and such that $\hat{y}_{\gamma_{k_l}}(t)\geq 0$, $\hat{w}_{\gamma_{k_l}}(t)=0$.

In what follows we restrict attention to problems where conditions  (\ref{pm1}) - (\ref{pm6}) for the approximating problems are such that   ${\rm meas\:}( E)=0$ or that 
$$
\frac{d}{dt}(-\partial_w\hat{g}_{\gamma_{k_l}}{r}+\gamma_{k_l}h_w(\hat{y}_{\gamma_{k_l}})_+r)=0,\;\; t\in D_E.
$$
The above derivative is understood in the following sense: for any $\phi$ 
$$ \frac{d}{dt} \phi(t^*)=\displaystyle \lim_{\stackrel{ t\;\longrightarrow\; t^*}{{\tiny D_E}} }\frac{\phi(t)- \phi(t^*)}{t-t^*}.$$


We now present  some examples to show how this main assumption can be verified.

\begin{example}\label{Example 1}

Consider the system
$$
\ddot{y}=u-N,\;\; N\geq 0,\;\; y\leq 0,\;\; u\in [-1,1].
$$
The corresponding approximating system has the form
\begin{eqnarray*}
&& \dot{y}=w,\\
&& \dot{w}=u-\gamma y_+w_+.
\end{eqnarray*}
We apply the necessary conditions \eqref{pm1}-\eqref{pm6} to the  the minimum time problem associated to the  system above (see problem  \eqref{0d}--\eqref{6d}). The adjoint system is
\begin{eqnarray*}
&& \dot{s}=\gamma h_yw_+r,\\
&& \dot{r}=-s + \gamma h_wy_+r.
\end{eqnarray*}
From the maximum condition \eqref{pm4} we deduce that
$$
u=\left\{
\begin{array}{cl}
 1, & r>0,\\
 -1, & r<0.
\end{array}
\right.
$$

Let $E$ be the set of all  Lebesgue point $t$  of 
$\dot{w}$,  $u$, $\dot{s}$, and $\dot{r}$ where  $y(t)\geq 0$ and  $w(t)=0$.
If $t\in D_E$, then $\dot{w}(t)=0$. Otherwise $t$ is an isolated point of $E$ or does not belong to $E$. Hence we have $u(t)=0$. The maximum condition implies that $r(t)=0$. Since $s^2+r^2\neq 0$, we obtain $s(t)\neq 0$. This implies that $\dot{r}(t)=-s(t)\neq 0$. Therefore $t\not\in D_E$. So, in this example ${\rm meas\:}(E)=0$.
\end{example}

\begin{example}\label{Example 2}

Consider a cylinder attached to a spring with a piston moving inside the  cylinder as in Figure \ref{fig:1}. Assume that we apply a  force $u\in[-1,1]$ to the piston. 
\begin{center}
\begin{figure}[h!]
\begin{center}
\pgfimage[width=9cm]{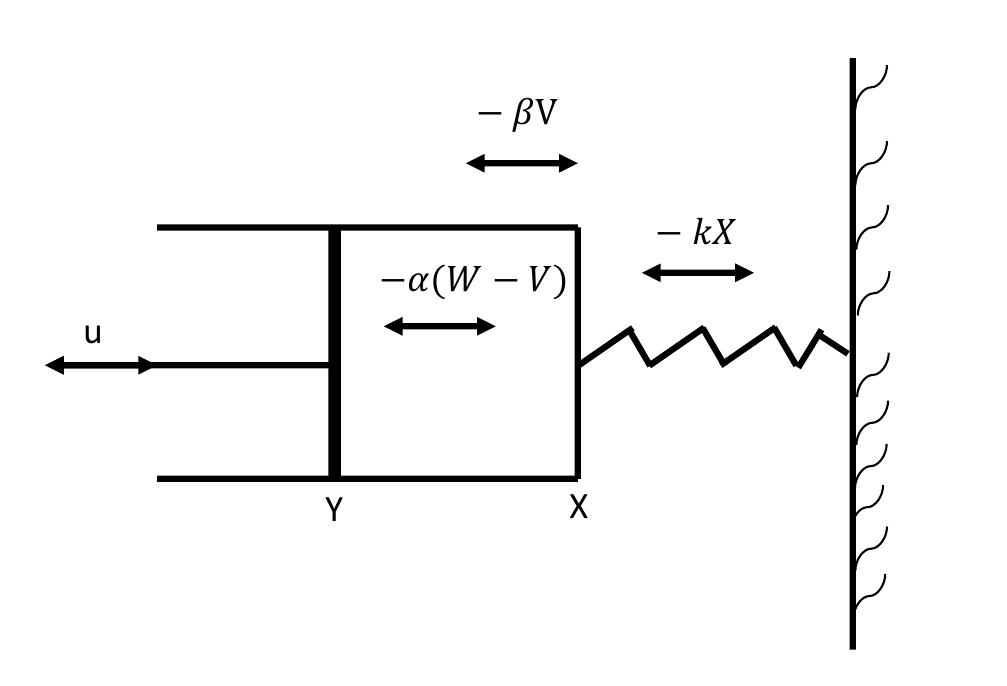}
\end{center}
\caption{\small{Here, $X$ is the position of the bottom of the cylinder and $Y$ is the position of the piston. The parameters $\alpha$ and $ \beta$ are the resistance coefficients while  $k$ is the stiffness of the spring. Moreover, we have $V=\dot X$ and $W=\dot Y$. }}
\label{fig:1}
\end{figure}
\end{center}

Informally, the equations of motion are 
\begin{eqnarray*}
&& M\ddot{X}=\alpha (\dot{Y}-\dot{X})-kX-\beta \dot X +N,\\
&& \ddot{Y}=\alpha (\dot{X}-\dot{Y}) + u- N,\\
&& N\geq 0,\;\;X\geq Y,\;\; u\in [-1,1],
\end{eqnarray*}
where $N$ is the reaction force.
In the above $M$ is the mass of the cylinder and  the mass of the piston is $1$.

Introducing new coordinates $x=MX+Y$ and $y=Y-X$ we obtain the system in canonical form
\begin{eqnarray*}
&& \dot{x}=v,\\
&& \dot{y}=w,\\
&& \dot{v}=-a (x-y)-c(v-w)+u,\\
&& \dot{w}=a (x-y)+c(v-w)-bw+u-2N,
\end{eqnarray*}
where
$$
a=\frac{k}{M+1},\;\; b=\alpha\frac{M+1}{M},\;\; c=\frac{\beta}{M+1}.
$$
The approximating system has the form
\begin{eqnarray*}
&& \dot{x}=v,\\
&& \dot{y}=w,\\
&& \dot{v}=-a (x-y)-c(v-w)+u,\\
&& \dot{w}=a (x-y)+c(v-w)-bw+u-\gamma y_+w_+,
\end{eqnarray*}
and the corresponding  adjoint system is
\begin{eqnarray}
&& \dot{q}=a(p-r),\label{ex2:1}\\
&& \dot{s}=-a(p-r)+\gamma h_yw_+r,\label{ex2:2}\\
&& \dot{p}=-q +c(p-r),\label{ex2:3}\\
&& \dot{r}=-s-c(p-r)+br + \gamma h_wy_+r.\label{ex2:4}
\end{eqnarray}
Moreover from the maximum condition we have
\begin{equation}\label{ex2:umax}
u=\left\{
\begin{array}{cl}
 1, & p+r>0,\\
 -1, & p+r<0.
\end{array}
\right.
\end{equation}
Let $E$ be the set of all Lebesgue points $t$  of $\dot{v}$, 
$\dot{w}$,  $u$, $\dot{q}$, $\dot{s}$, $\dot{p}$, and $\dot{r}$  where  $y(t)\geq 0$ and $w(t)=0$.

If $t\in D_E$, then $\dot{w}(t)=0$ and we have 
\begin{equation}
0=a(x(t)-y(t))+cv(t)+u(t).\label{e1}
\end{equation} 
Hence 
\begin{equation}
\dot{v}(t)=2u(t).\label{e2}
\end{equation}
 Consider a sequence $t_j\in D_E$, $j=1,2,\ldots$, converging to $t$. Since at all points $t_j$   equality 
 (\ref{e1}) holds, we get
 $$
 0=a\frac{x(t_j)-x(t)}{t_j-t}-a\frac{y(t_j)-y(t)}{t_j-t}+c\frac{v(t_j)-v(t)}{t_j-t}+\frac{u(t_j)-u(t)}{t_j-t}.
 $$
 Passing to the limit, we obtain
 $$
 0=av(t)+c\dot{v}(t)+\lim_{k\rightarrow\infty}\frac{u(t_j)-u(t)}{t_j-t}.
 $$
Seeking a contradiction, assume that  $p(t)+r(t)\neq 0$. It follows from \eqref{ex2:umax} that  $u(t)=u(t_j)=\pm 1$ and we get
 $$
 0=av(t)+c\dot{v}(t).
 $$
 In the same way we have
 $$
 0=av(t_j)+c\dot{v}(t_j).
 $$
From (\ref{e2}) we obtain
$$
0=av(t_j)\pm 2c.
$$
Passing to the limit in the  equality
$$
0=a\frac{v(t_j)-v(t)}{t_j-t},
$$
we get $0=\dot{v}(t)=2u(t)$. Hence by the maximum condition we have $p(t)+r(t)=0$. The same is true for the points $t_j$. Passing to the limit in the equality  
$$
0=\frac{p(t_j)-p(t)}{t_j-t}+\frac{r(t_j)-r(t)}{t_j-t},
$$
we conclude that $\dot{p}(t)+\dot{r}(t)=0$. In the same way we have $\dot{p}(t_j)+\dot{r}(t_j)=0$. Since $\ddot{p}(t)$ exists, we see from \eqref{ex2:3} that there exists $\ddot{r}(t)$ and 
$$
\ddot{r}(t)=-\ddot{p}(t)=\dot{q}(t)-c(\dot{p}(t)-\dot{r}(t))=-\dot{s}-c(\dot{p}(t)-\dot{r}(t)).
$$
On the other hand,  we have
$$
\ddot{r}(t)=-\dot{s}-c(\dot{p}(t)-\dot{r}(t))+\frac{d}{dt}(br+\gamma h_wy_+r).
$$
Hence 
$$
\frac{d}{dt}(br+\gamma h_wy_+r)=0,\;\; t\in D_E.
$$
This shows that  our main assumption holds for this example.
\end{example}

\section{Necessary conditions}

\begin{theorem}\label{main}

There exists a solution  $(\bar{x},\bar{y},\bar{v},\bar{w},\bar{u})$ to problem \eqref{0b}--\eqref{6b}  such that there exist 
$$
(q,\sigma, p,r,\mu)\in AC([0,\hat{T}],R^n)\times BV([0,\hat{T}],R)\times AC([0,\hat{T}],R^n)\times AC([0,\hat{T}],R)\times BV([0,\hat{T}],R)
$$
and a constant $\bar H$
 satisfying 

\begin{eqnarray}
&& \dot{q}=-(\nabla_x\bar{f})^*p-\nabla_x\bar{g}r,\label{pm1a}\\
&& d\sigma= -\langle \partial_y\bar{f},p\rangle dt-\partial_y\bar{g}rdt +d\mu,\label{pm2a}\\
&& \dot{p}=-q-(\nabla_v\bar{f})^*p-\nabla_v\bar{g}r,\label{pm3a}\\
&& \dot{r}=-\sigma-\langle \partial_w\bar{f},p\rangle-\partial_w\bar{g}r,\label{pm4a}\\
&& \max_{u\in U(t)}(\langle q, \bar{v}\rangle+\sigma \bar{w}+\langle p, 
f(\bar{x},\bar{y},\bar{v},\bar{w},u)\rangle
 +
rg(\bar{x},\bar{y},\bar{v},\bar{w},u))\nonumber\\
&& =\langle q, \bar{v}\rangle+\sigma \bar{w}+\langle p, 
f(\bar{x},\bar{y},\bar{v},\bar{w},\bar{u})\rangle
 +
rg(\bar{x},\bar{y},\bar{v},\bar{w},\bar{u})\equiv \bar H \geq 0, \label{pm5a}\\
&& |q(\hat{T})|^2+|\sigma(\hat{T})|^2+|p(\hat{T})|^2+|r(\hat{T})|^2=1. \label{pm6a}
\end{eqnarray}
where $\bar{f}=f(\bar{x},\bar{y},\bar{v},\bar{w},\bar{u})$ and $\bar{g}=g(\bar{x},\bar{y},\bar{v},\bar{w},\bar{u})$.
Moreover
$$
\bar{x}(\hat{T})=x_1,\;\; \bar{y}(\hat{T})=y_1,\;\;
\bar{v}(\hat{T})=v_1,\;\;
\bar{w}(\hat{T})=w_1,
$$
and
$d\mu=0$, whenever $y(t)<0$. 
\end{theorem}

\begin{proof}
Recall the definition \eqref{limit}.    For each   $k_l$, consider  the    approximating problem \eqref{0d}--\eqref{6d}.  The corresponding  necessary conditions for that   problem   are given by \eqref{pm1}--\eqref{pm6}. 

Fixing $k_l$, set 
$$
\xi=s+\langle \partial_w\hat{f}_{\gamma_{k_l}},p\rangle+\partial_w\hat{g}_{\gamma_{k_l}}r-{\gamma_{k_l}}h_w(\hat{y}_{\gamma_{k_l}})_+r.
$$
Under our  main assumption, we have 
\begin{eqnarray}
&& 
\dot{\xi}=
-\langle \partial_y\hat{f}_{\gamma_{k_l}},p\rangle-\partial_y\hat{g}_{\gamma_{k_l}}r+{\gamma_{k_l}}
h_y(\hat{w}_{\gamma_{k_l}})_+r
+\frac{d}{dt}\langle \partial_w\hat{f}_{\gamma_{k_l}},p\rangle +\frac{d}{dt}(\partial_w\hat{g}_{\gamma_{k_l}}r) \nonumber\\
&&
+
\left\{
\begin{array}{cl}
-\frac{d}{dt}(\partial_w\hat{g}_{\gamma_{k_l}}r), & t\in D_E,\\
0, & t\in\{ t~:~\hat{y}_{\gamma_{k_l}}(t)<0 \; {\rm or}\; \hat{w}_{\gamma_{k_l}}(t)<0\},\\
-{\gamma_{k_l}}h_w\hat{w}_{\gamma_{k_l}}r-{\gamma_{k_l}}
(\hat{y}_{\gamma_{k_l}})_+\dot{r}, & t\in\{ t~:~\hat{y}_{\gamma_{k_l}}(t)> 0 \; {\rm and}\; 
\hat{w}{\gamma_{k_l}}(t)>0\},
\end{array}
\right.\nonumber\\
&& 
=
-\langle \partial_y\hat{f}_{\gamma_{k_l}},p\rangle-\partial_y\hat{g}_{\gamma_{k_l}}r+
\frac{d}{dt}\langle \partial_w\hat{f}_{\gamma_{k_l}},p\rangle +\frac{d}{dt}(\partial_w\hat{g}_{\gamma_{k_l}}r) \nonumber\\
&&
+
\left\{
\begin{array}{cl}
-\frac{d}{dt}(\partial_w\hat{g}_{\gamma_{k_l}}r), & t\in D_E,\\
0, & t\in\{ t~:~y(t)<0 \; {\rm or}\; w(t)<0\},\\
-{\gamma_{k_l}}(\hat{y}_{\gamma_{k_l}})_+\dot{r}, & t\in\{ t~:~y(t)> 0 \; {\rm and}\; w(t)>0\},
\end{array}
\right. \label{xi}
\end{eqnarray}
and
$$
\dot{r}=-\xi.
$$
(In the above, notice that   ${\rm meas} (\{ t~:~\hat{y}_{\gamma_{k_l}}(t)=0,\; \hat{w}_{\gamma_{k_l}}(t)\neq 0\})=0$.) From these equations and 
recalling that $\langle a,b\rangle \geq -|\langle a,b\rangle| \geq -\frac{|a|^2+|b|^2}{2}$   for any two vectors $a$ and $b$,  we get
$$
\frac{1}{2}\frac{d}{dt}\xi^2\geq -M_5(|q|^2+\xi^2+|p|^2+r^2).
$$
Using the above inequality and multiplying  \eqref{0d} by $q$, \eqref{2d} by $p$ and  $\dot r=-\xi$  by $r$,  we deduce that  
$$
\frac{1}{2}\frac{d}{dt}(|q|^2+\xi^2+|p|^2+r^2)\geq -M_6(|q|^2+\xi^2+|p|^2+r^2).
$$
Appealing now to  the Gronwall inequality and the nontriviality condition
$$
|q(\hat{T}_l)|^2+|s(\hat{T}_l)|^2+|p(\hat{T}_l)|^2+|r(\hat{T}_l)|^2=1,
$$
  we see that $|q|^2+\xi^2+|p|^2+r^2$ is bounded
$$
|q|^2+\xi^2+|p|^2+r^2\leq M_7
$$
and this estimate does not depend on $\gamma$.

Multiplying (\ref{xi}) by ${\rm sign\:}(\xi(t))$ and integrating, we get
\begin{equation}\label{xi:l}
|\xi(\hat{T}_l)|-|\xi(0)|\geq - M_8 +\int_0^{\hat{T}_l}{\gamma_{k_l}}(\hat{y}_{\gamma_{k_l}})_+|\xi(t)|dt.
\end{equation}

Up to now we considered $k_l$ fixed and so the adjoint variable had no subscript. From now on, we want to focus on the case where $l\to \infty$ and so we need to identify the adjoint variables with the subscript $\gamma_{k_l}$.  Set 
$$
\sigma_{\gamma_{k_l}}=s_{\gamma_{k_l}}-
{\gamma_{k_l}}h_w(\hat{y}_{\gamma_{k_l}})_+r_{\gamma_{k_l}}
$$
and 
$$\dot \mu_{\gamma_{k_l}}dt = 
\left\{
\begin{array}{cl}
-\frac{d}{dt}(\partial_w\hat{g}_{\gamma_{k_l}}r_{\gamma_{k_l}})dt, & t\in (D_E)_{\gamma_{k_l}},\\
0, & t\in\{ t~:~\hat{y}_{\gamma_{k_l}}(t)<0 \; {\rm or}\; \hat{w}_{\gamma_{k_l}}(t)<0\},\\
{\gamma_{k_l}}(\hat{y}_{\gamma_{k_l}})_+\xi_{\gamma_{k_l}} dt, & t\in\{ t~:~\hat{y}_{\gamma_{k_l}}(t)> 0 \; {\rm and}\; \hat{w}_{\gamma_{k_l}}(t)>0\}.
\end{array}
\right.
$$
Recall that $\sigma_{\gamma_{k_l}}=s_{\gamma_{k_l}}$ whenever $w_{\gamma_{k_l}}<0$. From \eqref{pm1}--\eqref{pm6} we have

\begin{eqnarray}
&& \dot q_{\gamma_{k_l}}=-(\nabla_x\hat{f}_{\gamma_{k_l}})^*p_{\gamma_{k_l}}-
\nabla_x\hat{g}_{\gamma_{k_l}}r_{\gamma_{k_l}},\label{pmm1}\\
&&\dot\sigma_ {\gamma_{k_l}}=-\langle \partial_y\hat{f}_{\gamma_{k_l}},p_{\gamma_{k_l}}\rangle-
\partial_y\hat{g}_{\gamma_{k_l}}r_{\gamma_{k_l}}+\dot \mu_{\gamma_{k_l}}\label{pmm2}\\
&& \dot p_{\gamma_{k_l}}=-q_{\gamma_{k_l}}-(\nabla_v\hat{f}_{\gamma_{k_l}})^*p_{\gamma_{k_l}} -\nabla_v\hat{g}_{\gamma_{k_l}}r_{\gamma_{k_l}} ,\label{pmm3}\\
&& \dot r_{\gamma_{k_l}}= \sigma_{\gamma_{k_l}}-\langle \partial_w\hat{f}_{\gamma_{k_l}},p_{\gamma_{k_l}} \rangle-\partial_w\hat{g}_vr_{\gamma_{k_l}} ,\label{pmm4}\\
&& \max_{u\in U}(\langle q_{\gamma_{k_l}} , \hat{v}_{\gamma_{k_l}}\rangle+\sigma_ {\gamma_{k_l}}\hat{w}_{\gamma_{k_l}}+\langle p_{\gamma_{k_l}} , 
f(\hat{x}_{\gamma_{k_l}},\hat{y}_{\gamma_{k_l}},\hat{v}_{\gamma_{k_l}},\hat{w}_{\gamma_{k_l}},u)\rangle \nonumber\\
&& \nonumber  \qquad   \qquad +
r_{\gamma_{k_l}}  g(\hat{x}_{\gamma_{k_l}},\hat{y}_{\gamma_{k_l}},\hat{v}_{\gamma_{k_l}},\hat{w}_{\gamma_{k_l}},u))\nonumber\\
&& =\langle q_{\gamma_{k_l}} , \hat{v}_{\gamma_{k_l}}\rangle+\sigma_ {\gamma_{k_l}}\hat{w}_{\gamma_{k_l}}+
\langle p_{\gamma_{k_l}} , \hat{f}_{\gamma_{k_l}} \rangle +
r_{\gamma_{k_l}} \hat{g}_{\gamma_{k_l}}\equiv \hat H_l \geq 0, \label{pmm5}\\[2mm]
&& |q_{\gamma_{k_l}} (\hat{T}_l)|^2+|\sigma_{\gamma_{k_l}} (\hat{T}_l)|^2+|p_{\gamma_{k_l}} (\hat{T}_l)|^2+|r_{\gamma_{k_l}} (\hat{T}_l)|^2=1.\label{pmm6}
\end{eqnarray}

 Consider the interval $[0,\hat T_1]$. From \eqref{xi:l} we see that, without loss of generality, the measures $ \dot \mu_{\gamma_{k_l}}dt $
converge in weak* topology of  to some measure  $d\mu^1$,  $\mu^1\in BV([0,\hat{T}_1],R)$. Hence (see \eqref{xi}), without loss of generality, $\xi_{\gamma_{k_l}}$ also
converges in weak* topology of $BV([0,\hat{T}_1],R)$ to some function $\xi^1(t)$. Therefore the sequence
$
\sigma_{\gamma_{k_l}}$
also converges in weak* topology of $BV([0,\hat{T}_1],R)$ to some function $\sigma^1(t)$ (we do not relabel).   The functions $(q_{\gamma_{k_l}},\sigma_{\gamma_{k_l}},p_{\gamma_{k_l}},r_{\gamma_{k_l}})$ are bounded in $[0,\hat{T}_1]$. It follows from \eqref{pmm1}, \eqref{pmm3} and \eqref{pmm4} that $(\dot q_{\gamma_{k_l}},\dot p_{\gamma_{k_l}},\dot r_{\gamma_{k_l}})$
are also bounded. Appealing to the Arzela-Ascoli Theorem we conclude that, without loss of generality, $(q_{\gamma_{k_l}},p_{\gamma_{k_l}},r_{\gamma_{k_l}})$ uniformly converges to some function $(q,p,r)$ in $[0,\hat{T}_1]$ and the sequence $\sigma_{\gamma_{k_l}}$ converges pointwisely to $\sigma$ by Helly's first Theorem. Integrating \eqref{pmm1}, \eqref{pmm3} and \eqref{pmm4} over $[t,\hat{T}_1]$, for $t<\hat T_1$ we have
\begin{eqnarray}
&&  q_{\gamma_{k_l}}(t)=q_{\gamma_{k_l}}(\hat T_1)-\int_t^{\hat T_1}(-(\nabla_x\hat{f}_{\gamma_{k_l}})^*p_{\gamma_{k_l}}-
\nabla_x\hat{g}_{\gamma_{k_l}}r_{\gamma_{k_l}})ds,\label{int1}\\
&& p_{\gamma_{k_l}}(t)=p_{\gamma_{k_l}}(\hat T_1)-\int_t^{\hat T_1}
(-q_{\gamma_{k_l}}-(\nabla_v\hat{f}_{\gamma_{k_l}})^*p_{\gamma_{k_l}} -\nabla_v\hat{g}_{\gamma_{k_l}}r_{\gamma_{k_l}})ds ,\label{int3}\\
&&  r_{\gamma_{k_l}}(t)= r_{\gamma_{k_l}}(\hat T_1)-\int_t^{\hat T_1}(-\sigma_{\gamma_{k_l}}-\langle \partial_w\hat{f}_{\gamma_{k_l}},p_{\gamma_{k_l}} \rangle-\partial_w\hat{g}_vr_{\gamma_{k_l}})ds \label{int4}
\end{eqnarray}
It is easy to conclude that 
$(x_{\gamma_{k_l}},y_{\gamma_{k_l}},v_{\gamma_{k_l}},w_{\gamma_{k_l}},u_{\gamma_{k_l}})$
converge to some 
$(\bar x, \bar y, \bar v, \bar w, \bar u)$, where   $\bar x \in AC([0,\hat T_1],R^n)$,   $\bar y  \in AC([0,\hat T_1],R)$,  $\bar v  \in AC([0,\hat T_1],R^n)$, $\bar w  \in  BV([0,\hat T_1],R)$ and $ \bar \mu \in  L_\infty(([0,\hat T_1],R^m)$,
and, moreover,  $(x_{\gamma_{k_l}},y_{\gamma_{k_l}},v_{\gamma_{k_l}})$ converge uniformly,  $w_{\gamma_{k_l}}$ converges in the weak* topology in $BV$ and $u_{\gamma_{k_l}}$ converges in the weak* topology in $L_\infty$.
Passing to the limit in (\ref{int1}) - (\ref{int4}) we get
\begin{eqnarray}
&& q(t)=q(\hat T_1) -\int_t^{\hat T_1}(-(\nabla_x\bar{f})^*p-\nabla_x\bar{g}r) ds ,\label{int1a}\\
&& p(t)=p(\hat T_1) -\int_t^{\hat T_1}(-q-(\nabla_v\bar{f})^*p-\nabla_v\bar{g}r)ds,\label{int3a}\\
&& r(t)=r(\hat T_1) -\int_t^{\hat T_1}(-\sigma-\langle \partial_w\bar{f},p\rangle-\partial_w\bar{g}r)ds,\label{int4a}
\end{eqnarray}
where $\bar{f}=f(\bar{x},\bar{y},\bar{v},\bar{w},\bar{u})$ and $\bar{g}=g(\bar{x},\bar{y},\bar{v},\bar{w},\bar{u})$. It follows that $(q,p,r)$ are absolutely continuous.

Consider now any continuous function $\sigma^*$ defined in $[0,\hat T_1]$. Then we have
$$\int_0^{\hat T_1} \sigma^* \dot\sigma_ {\gamma_{k_l}}ds=\int_0^{\hat T_1} \sigma^*(-\langle \partial_y\hat{f}_{\gamma_{k_l}},p_{\gamma_{k_l}}\rangle-
\partial_y\hat{g}_{\gamma_{k_l}}r_{\gamma_{k_l}}+\dot \mu_{\gamma_{k_l}})ds.$$
Passing to the limit in above equality we get \eqref{pm2a} in the interval $[0,\hat T_1]$.

Next, we consider the interval $[0,\hat T_2]$. From the previous subsequences involving the control processes and dual variables we extract another sequence where the convergence holds in the interval $[0,\hat T_2]$, with $\hat T_2 >\hat T_1$. We repeat  this process for all intervals $[0,\hat T_3]$, $[0,\hat T_4]$,  etc. Then, we extract the diagonal subsequence to get convergence in the interval $[0, \hat T]$. 

We now concentrate on \eqref{pmm5}. Since $f$ and $g$ are affine with respect to $u$,  we have 
$$
\begin{array}{c}
\langle q_{\gamma_{k_l}} , \hat{v}_{\gamma_{k_l}}\rangle+\sigma_ {\gamma_{k_l}}\hat{w}_{\gamma_{k_l}}+\langle p_{\gamma_{k_l}} , 
f(\hat{x}_{\gamma_{k_l}},\hat{y}_{\gamma_{k_l}},\hat{v}_{\gamma_{k_l}},\hat{w}_{\gamma_{k_l}},u)\rangle+
r_{\gamma_{k_l}}  g(\hat{x}_{\gamma_{k_l}},\hat{y}_{\gamma_{k_l}},\hat{v}_{\gamma_{k_l}},\hat{w}_{\gamma_{k_l}},u)\\[2mm]
=H^0_l+\langle h_l,u\rangle
\end{array}
$$
and 
$$
\begin{array}{c}
\langle q_{\gamma_{k_l}} , \hat{v}_{\gamma_{k_l}}\rangle+\sigma_ {\gamma_{k_l}}\hat{w}_{\gamma_{k_l}}+
\langle p_{\gamma_{k_l}} , \hat{f}_{\gamma_{k_l}} \rangle +
r_{\gamma_{k_l}} \hat{g}_{\gamma_{k_l}}\equiv H_l \\[2mm]
= H^0_l+\langle  h_l,\hat u_l\rangle\geq 0,
\end{array}
$$
for some   $H^0_l$ and $h_l$. Thus, \eqref{pmm5} reads
\begin{equation}\label{Ham:l}
\max_{u \in U}\left\{H^0_l+\langle  h_l,u\rangle \right\}=H^0_l+\langle  h_l,\hat u_l\rangle.
\end{equation}
We now take limits in \eqref{Ham:l}, i.e., we want to prove that   \eqref{Ham:l} converges  to
\begin{equation}\label{Ham:l2}\max_{u \in U}\left\{\bar H^0+\langle  \bar h,u\rangle \right\}=\bar H^0+\langle  \bar h,\bar u\rangle.
\end{equation}

It is an easy matter to see that $H^0_l$ and $h_l$ are bounded functions converging pointwisely to some functions $\bar H^0 $ and $\bar h $.
Without loss of generality, we deduce that $\hat H_l \to \bar  H$. Next, we show that 
$\bar H^0+\langle  \bar h,\bar u \rangle=\bar H$ almost everywhere. Seeking a contradiction, assume that 
$$
\text{meas}(I^>) =\text{meas}\left\{ t:~\bar H^0+\langle  \bar h ,\bar u\rangle> \bar H \right\}>0.
$$
Integrating, we get 
$$\begin{array}{c}\bar H \text{meas}(I^>) = \displaystyle\lim_{l\to \infty}  \displaystyle \int_{I^>} \left(H^0_l+\langle  h_l,\hat u_l\rangle\right)  dt = \displaystyle \lim_{l\to \infty}  \displaystyle\int_{I^>} \left(H^0_l+\langle  h_l-\bar h,\hat u_l\rangle +  \langle \bar h,\hat u_l\rangle\right) dt \\[2mm]
= \displaystyle\int_{I^>} \left(\bar H^0+  \langle \bar h,\bar u\rangle\right) dt ~> ~\bar H \text{meas}(I^>).
\end{array}
$$
This means that $\bar H^0+\langle  \bar h,\bar u \rangle\leq \bar H$ almost everywhere. Analogously we show that $\bar H^0+\langle  \bar h,\bar u \rangle\geq \bar H$, proving our claim that $\bar H^0+\langle  \bar h,\bar u \rangle= \bar H$.

We now assume that there exists an interval $\tilde I \subset [0, \hat T]$ of positive measure  such that 
$\max_{u \in U}\left\{\bar H^0+\langle  \bar h,u\rangle \right\} > \bar H^0+\langle  \bar h,\bar u\rangle$ for all $t\in \tilde I$.  Then, there exists an admissible control $\tilde u$ such that
$$\bar H^0+\langle  \bar h,\tilde u \rangle > \bar H^0+\langle  \bar h,\bar u\rangle=\bar H~ \text{ for all } t\in \tilde I.$$
It follows that there exists $\epsilon >0$ such that
$$\begin{array}{c}
\bar H \text{meas}(\tilde I) +2\epsilon < \displaystyle \int_{\tilde I} \left(\bar H^0+\langle  \bar h,\tilde u \rangle\right) dt\\[2mm]
= \displaystyle \int_{\tilde I} \left( H^0_l+\langle  h_l,\tilde u \rangle\right) dt +\displaystyle \int_{\tilde I} \left(\bar H^0-H^0_l+\langle  \bar h-h_l,\tilde u \rangle\right) dt\\[2mm]
 \leq  H_l~\text{meas}(\tilde I)+\epsilon \leq \bar H \text{meas}(\tilde I) +2\epsilon.
\end{array}
$$
So we deduce that measure of $\tilde I$ is $0$, proving \eqref{pm5a}. The theorem is then proved.
\end{proof}

The main disadvantage of these necessary conditions is that, in general, we cannot guarantee the non triviality condition 
\begin{equation}
|q(t)|^2+|\sigma(t)|^2+|p(t)|^2+|r(t)|^2>0 \label{pmN}
\end{equation}
for all $t\in [0,\hat{T}]$. Below we present a special case where (\ref{pmN}) is satisfied for all $t\in [0,\hat{T}]$.

\section{The case when the target set is a stable neighborhood of an equilibrium point}

To better explain the main idea we consider a linear control system
\begin{equation}
\dot{z}=Az+Bu,\;\; u\in U,\label{ls1}
\end{equation}
where $U\subset R^m$ is a convex compact set containing zero in its interior. Assume that the system is controllable. Then there exists a 
linear feedback $u=Cz$ such that the zero equilibrium position of the system
$$
\dot{z}=(A+BC)z
$$
is asymptotically stable. Then there exists a positive definite matrix $V$ such that 
$$
\langle Vz,(A+BC)z\rangle\leq-\frac{1}{2}|z|^2.
$$
Now suppose that $z(t)$ is a time-optimal trajectory of system (\ref{ls1}) in a problem with the terminal set 
$S=\{ z~:~\langle z,Vz\rangle\leq\epsilon\}$. (If $\epsilon>0$ is sufficiently small, then $Cz\in U$, $z\in S$.) Then the transversality condition at the final point $z(T)\in S$ is $p(T)=-\frac{Vz(T)}{|Vz(T)|}$. Let $u$ be the optimal control. By the maximum principle we have
$$
\langle p(T), Az(T)+Bu(T)\rangle\geq \langle p(T), (A+BC)z(T)\rangle\geq 
\frac{|z(T)|^2}{2|Vz(T)|}\geq \rho(\epsilon)>0.
$$
This guarantees that the Hamiltonian, which is constant along the optimal trajectory, is strictly positive. For brevity we shall call such a set 
$S$ {\em a stable neighborhood of an equilibrium point}.

The same is true for  nonlinear systems where the terminal point $(x_1,y_1,0,0)$  is an equilibrium position and  the linearization at this point is controllable.

\vspace{5mm}


Let  $(x_1,y_1,0,0)$, $y_1<0$, be an equilibrium point  of  control system \eqref{1b}-\eqref{4b} and 
$$
S=\{ (x,y,v,w)~:~\langle (x-x_1,y-y_1,v,w),V(x-x_1,y-y_1,v,w)\rangle\leq\epsilon\}
$$
 be a small stable neighborhood of this point.
Let $(\tilde{x},\tilde{y},\tilde{v},\tilde{w},\tilde{u})$ be a global optimal solution of the problem
\begin{eqnarray}
&& T\rightarrow\min,\label{0cc}\\ 
&& \dot{x}=v,\label{1cc}\\
&& \dot{y}=w,\label{2cc}\\
&& \dot{v}=f(x,y,v,w,u),\label{3cc}\\
&& dw=g(x,y,v,w,u)dt-d\nu,\;\; d\nu\geq 0,\;\; d\nu=0,\; y<0,\label{4cc}\\
&& x(0)=x_0,\;y(0)=y_0<0,\;v(0)=v_0,\;w(0)=w_0,\label{5cc}\\
&& (x(T),y(T),v(T),w(T))\in S.\label{6cc}
\end{eqnarray} 
Let $\hat{T}$ be the optimal time. 
There exist a sequence ${\gamma_{k_l}}\rightarrow\infty$ and a sequence $u_{\gamma_{k_l}}$ converging to $\tilde{u}$ in the weak* topology of 
$L_{\infty}([0,\hat{T}],R^m)$ such that 
$$
(\tilde{x},\tilde{y},\tilde{v},\tilde{w})=\lim_{l\rightarrow\infty}(x_{\gamma_{k_l}},y_{\gamma_{k_l}},v_{\gamma_{k_l}},w_{\gamma_{k_l}}),
$$
where $(x_{\gamma_{k_l}},y_{\gamma_{k_l}},v_{\gamma_{k_l}},w_{\gamma_{k_l}})$ is the solution to the Cauchy problem
\begin{eqnarray}
&& \dot{x}_{\gamma_{k_l}}=v_{\gamma_{k_l}},\label{1ccc}\\
&& \dot{y}_{\gamma_{k_l}}=w_{\gamma_{k_l}},\label{2ccc}\\
&& \dot{v}_{\gamma_{k_l}}=f(x_{\gamma_{k_l}},y_{\gamma_{k_l}},v_{\gamma_{k_l}},w_{\gamma_{k_l}},u_{\gamma_{k_l}}),\label{3ccc}\\
&& \dot{w}_{\gamma_{k_l}}=g(x_{\gamma_{k_l}},
y_{\gamma_{k_l}},v_{\gamma_{k_l}},w_{\gamma_{k_l}},u_{\gamma_{k_l}})-{\gamma_{k_l}} (y_{\gamma_{k_l}})_+(w_{\gamma_{k_l}})_+,\label{4ccc}\\
&& x_{\gamma_{k_l}}(0)=x_0,\;y_{\gamma_{k_l}}(0)=y_0,\;v_{\gamma_{k_l}}(0)=v_0,\;w_{\gamma_{k_l}}(0)=w_0.\label{5ccc}
\end{eqnarray}
We have 
$$
(x_{\gamma_{k_l}},y_{\gamma_{k_l}},v_{\gamma_{k_l}},w_{\gamma_{k_l}})(T_l)\in (x_1,y_1,0,0)+(1+\delta_l)(S-(x_1,y_1,0,0)), 
$$
 where $T_l\uparrow \hat{T}$ and 
$\delta_l\downarrow 0$. 
Consider the following time-optimal control problem:
\begin{eqnarray}
&& T\rightarrow\min,\label{0dd}\\
&& \dot{x}=v,\label{1dd}\\
&& \dot{y}=w,\label{2dd}\\
&& \dot{v}=f(x,y,v,w,u),\label{3dd}\\
&& \dot{w}=g(x,y,v,w,u)-{\gamma_{k_l}}y_+w_+,\label{4dd}\\
&& x(0)=x_0,\;y(0)=y_0,\;v(0)=v_0,\;w(0)=w_0,\label{5dd}\\
&&  (x(T),y(T),v(T),w(T))\in (x_1,y_1,0,0)+(1+\delta_l)(S-(x_1,y_1,0,0)).\label{6dd}
\end{eqnarray} 
 This problem has an admissible process $(x_{\gamma_{k_l}},y_{\gamma_{k_l}},v_{\gamma_{k_l}},w_{\gamma_{k_l}}, u_{\gamma_{k_l}})$.
Denote by 
$(\hat{x}_{\gamma_{k_l}},\hat{y}_{\gamma_{k_l}},\hat{v}_{\gamma_{k_l}},\\ \hat{w}_{\gamma_{k_l}},\hat{u}_{\gamma_{k_l}})$ an optimal solution to this problem. The corresponding optimal time is $\hat{T}_l\leq T_l$.

\vspace{5mm}

Applying necessary conditions of optimality to this problem we see that there 
exist $(q,s,p,r)\in AC([0,\hat{T}],R^{n+1+n+1})$ and a constant $H_l$  satisfying the following conditions:
\begin{eqnarray}
&& \dot{q}=-(\nabla_x\hat{f}_{\gamma_{k_l}})^*p-\nabla_x\hat{g}_{\gamma_{k_l}}r,\label{pm1e}\\
&&\dot{s}=-\langle \partial_y\hat{f}_{\gamma_{k_l}},p\rangle-\partial_y\hat{g}_{\gamma_{k_l}}r+{\gamma_{k_l}}h_y(\hat{w}_{\gamma_{k_l}})_+r,\label{pm2e}\\
&& \dot{p}=-q-(\nabla_v\hat{f}_{\gamma_{k_l}})^*p-\nabla_v\hat{g}_{\gamma_{k_l}}r,\label{pm3e}\\
&& \dot{r}=-s-\langle \partial_w\hat{f}_{\gamma_{k_l}},p\rangle-\partial_w\hat{g}_{\gamma_{k_l}}r+{\gamma_{k_l}}h_w(\hat{y}_{\gamma_{k_l}})_+r,\label{pm4e}\\
&& \max_{u\in U}(\langle q, \hat{v}_{\gamma_{k_l}}\rangle+s\hat{w}_{\gamma_{k_l}}+\langle p, 
f(\hat{x}_{\gamma_{k_l}},\hat{y}_{\gamma_{k_l}},\hat{v}_{\gamma_{k_l}},\hat{w}_{\gamma_{k_l}},u)\rangle\\
&& \nonumber +
r(g(\hat{x}_{\gamma_{k_l}},\hat{y}_{\gamma_{k_l}},\hat{v}_{\gamma_{k_l}},\hat{w}_{\gamma_{k_l}},u)-{\gamma_{k_l}}(\hat{y}_{\gamma_{k_l}})_+(\hat{w}_{\gamma_{k_l}})_+))\nonumber\\
&& =\langle q, \hat{v}_{\gamma_{k_l}}\rangle+s\hat{w}_{\gamma_{k_l}}+\langle p, 
\hat{f}_{\gamma_{k_l}}\rangle +
r(\hat{g}_-{\gamma_{k_l}}(\hat{y}_{\gamma_{k_l}})_+(\hat{w}_{\gamma_{k_l}})_+)\equiv H_l\geq\rho(\epsilon)> 0, \label{pm5e}\\
&& (q,s,p,r)(\hat{T}_l)=-\frac{V(\hat{x}_{\gamma_{k_l}}-x_1,\hat{y}_{\gamma_{k_l}}-y_1,\hat{v}_{\gamma_{k_l}},\hat{w}_{\gamma_{k_l}})}{|V(\hat{x}_{\gamma_{k_l}}-x_1,\hat{y}_{\gamma_{k_l}}-y_1,\hat{v}_{\gamma_{k_l}},\hat{w}_{\gamma_{k_l}})|}~,\label{pm6e}
\end{eqnarray}
where $\hat{f}_{\gamma_{k_l}}=f(\hat{x}_{\gamma_{k_l}},\hat{y}_{\gamma_{k_l}},\hat{v}_{\gamma_{k_l}},\hat{w}_{\gamma_{k_l}},\hat{u}_{\gamma_{k_l}})$, 
$\hat{g}_{\gamma_{k_l}}=g(\hat{x}_{\gamma_{k_l}},\hat{y}_{\gamma_{k_l}},\hat{v}_{\gamma_{k_l}},\hat{w}_{\gamma_{k_l}},\hat{u}_{\gamma_{k_l}})$.
The functions
$$
h_y(t)=\left\{
\begin{array}{cl}
0, & \hat{y}_{\gamma_{k_l}}(t)<0,\\
h\in [0,1], & \hat{y}_{\gamma_{k_l}}(t)=0,\\
1, & \hat{y}_{\gamma_{k_l}}(t)>0,
\end{array}
\right.
$$
and
$$
h_w(t)=\left\{
\begin{array}{cl}
0, & \hat{w}_{\gamma_{k_l}}(t)<0,\\
h\in [0,1], & \hat{w}_{\gamma_{k_l}}(t)=0,\\
1, & \hat{w}_{\gamma_{k_l}}(t)>0,
\end{array}
\right.
$$
are measurable.

Arguing as above we deduce that the following result holds.

\begin{theorem}\label{theo:examp}
There exists a time-optimal control process $(\bar{x},\bar{y},\bar{v},\bar{w},\bar{u})$,  
$$
(q,\sigma, p,r,\mu)\in AC([0,\hat{T}],R^n)\times BV([0,\hat{T}],R)\times AC([0,\hat{T}],R^n)\times AC([0,\hat{T}],R)\times BV([0,\hat{T}],R)
$$ and a constant $\bar H$ 
 satisfying 

\begin{eqnarray}
&& \dot{q}=-(\nabla_x\bar{f})^*p-\nabla_x\bar{g}r,\label{pm1aa}\\
&& d\sigma= -\langle \partial_y\bar{f},p\rangle dt-\partial_y\bar{g}rdt +d\mu,\label{pm2aa}\\
&& \dot{p}=-q-(\nabla_v\bar{f})^*p-\nabla_v\bar{g}r,\label{pm3aa}\\
&& \dot{r}=-\sigma-\langle \partial_w\bar{f},p\rangle-\partial_w\bar{g}r,\label{pm4aa}\\
&& \max_{u\in U}(\langle q, \bar{v}\rangle+\sigma \bar{w}+\langle p, 
f(\bar{x},\bar{y},\bar{v},\bar{w},u)\rangle
 +
rg(\bar{x},\bar{y},\bar{v},\bar{w},u))\nonumber\\
&& =\langle q, \bar{v}\rangle+\sigma \bar{w}+\langle p, 
f(\bar{x},\bar{y},\bar{v},\bar{w},\bar{u})\rangle
 +
rg(\bar{x},\bar{y},\bar{v},\bar{w},\bar{u})\equiv \bar H \geq\rho(\epsilon)> 0, \label{pm5aa}\\
&& |q(\hat{T})|^2+|\sigma(\hat{T})|^2+|p(\hat{T})|^2+|r(\hat{T})|^2=1. \label{pm6aa}
\end{eqnarray}
where $\bar{f}=f(\bar{x},\bar{y},\bar{v},\bar{w},\bar{u})$ and $\bar{g}=g(\bar{x},\bar{y},\bar{v},\bar{w},\bar{u})$.
Moreover
$$
(\bar{x},\bar{y},\bar{v},\bar{w},\bar{u})(\hat{T})\in S,
$$
and
$d\mu=0$, whenever $y(t)<0$. 
\end{theorem}

\section{Examples}

We now go back to the systems introduced  in the Examples \ref{Example 1} and \ref{Example 2}.

\subsection{Example 1}

\begin{center}
\begin{figure}[h!]
\begin{center}
\pgfimage[width=8cm]{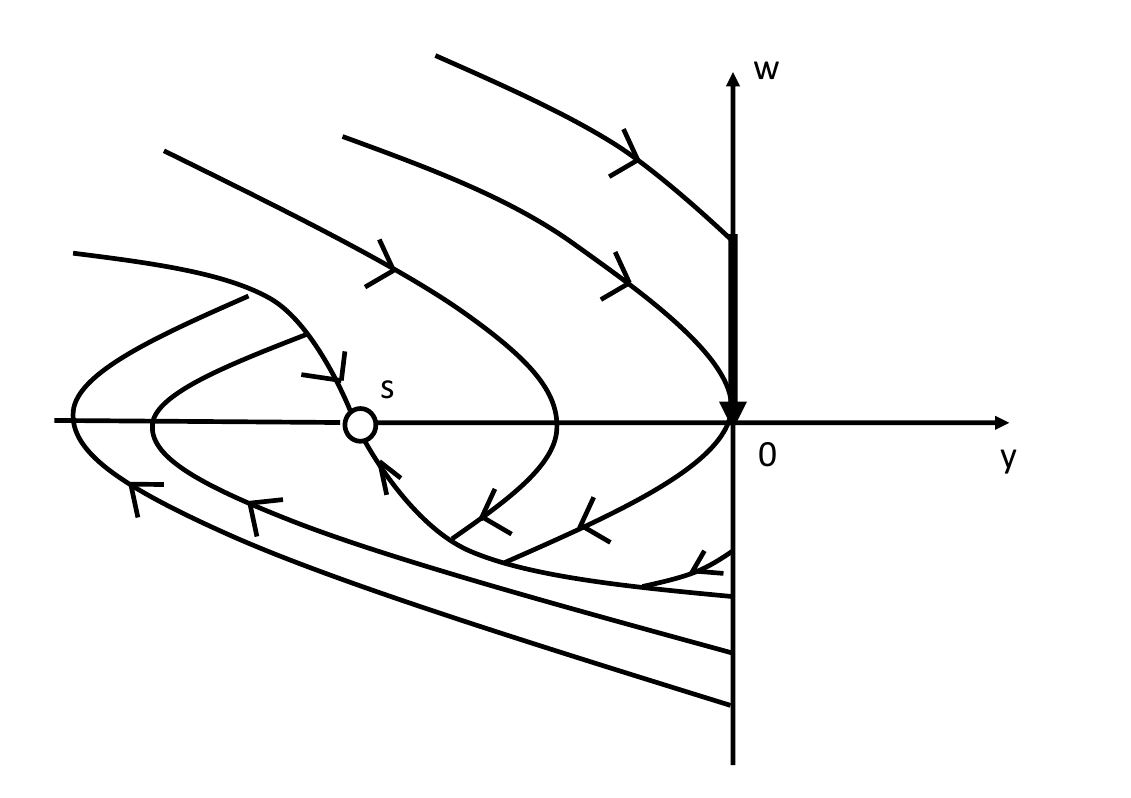}
\end{center}
\caption{\small{Time-optimal trajectories for  Example 1. The bold vertical arrow along the $w$ axis shows the jump of velocity when the trajectory reaches the boundary of the set $C=\left\{(y,w):~y\leq 0\right\}$.}}
\label{fig:2}
\end{figure}
\end{center}

Consider the problem associated to the system in Example  \ref{Example 1}:
\begin{eqnarray*}
&& T\rightarrow\min,\\
&& \ddot{y}=u-N,\;\; N\geq 0,\;\; y\leq 0,\;\; u\in [-1,1],\\
&& y(0)=y_0<0,\; \dot{y}(0)=v_0,\; (y(T),\dot{y}(T))\in S,
\end{eqnarray*}
where $S$ is a small stable neighborhood of a point $(y_1,0)$, $y_1<0$. 
Applying the necessary conditions we get
\begin{eqnarray*}
&& d\sigma=d\mu,\\
&& \dot{r}=-\sigma.
\end{eqnarray*}
Moreover from the maximum condition we have
$$
u=\left\{
\begin{array}{cl}
 1, & r>0,\\
 -1, & r<0.
\end{array}
\right.
$$
If $y<0$ we have $d\sigma=0$. Hence $r=at+b$, $a^2+b^2>0$. Therefore the control is constant or  changes sign just ones. The motion along the boundary $y=0$ is the passage from a positive velocity to zero. The optimal trajectories are show in the Figure \ref{fig:2}.

\subsection{Example 2}

Consider again a piston moving inside a cylinder as in Example \ref{Example 2} . 
Equations of motion are 
\begin{eqnarray*}
&& M\ddot{X}=\alpha (\dot{Y}-\dot{X})-kX-\beta \dot X+N,\\
&& \ddot{Y}=\alpha (\dot{X}-\dot{Y}) + u- N,\\
&& N\geq 0,\;\;X\geq Y,\;\; u\in [-1,1].
\end{eqnarray*}
We introduce the coordinates $x=MX+Y$ and $y=Y-X$ and consider the problem
\begin{eqnarray*}
&& T\rightarrow\min,\\
&& \dot{x}=v,\\
&& \dot{y}=w,\\
&& \dot{v}=-a (x-y)-c(v-w)+u,\\
&& \dot{w}=a (x-y)+c(v-w)-bw+u-N,\;\; N\geq 0,\;\; y\leq 0,\;\; u\in [-1,1],\\
&& x(0)=x_0,\; y(0)=y_0<0,\; v(0)=v_0<0,\; w(0)=w_0>0,\\
&& (x,y,v,w)(T)\in S,
\end{eqnarray*}
where $S$ is a small stable neighborhood of a point $(y_1,y_1,0,0)$, $y_1<0$. (The system is controllable near this point.)
Here $x_0=MX_0+Y_0$, $y_0=Y_0-X_0$, $X_0>0$, $Y_0<0$, $v_0=M\dot{X}_0<0$,  $w_0=-\dot{X}_0$, and 
$y_1=Y_0$. 
Applying the necessary conditions we see that there exist nontrivial $(q,\sigma,p,r)$ and $\mu$ satisfying
\begin{eqnarray*}
&& \dot{q}=a(p-r),\\
&& d\sigma=-a(p-r)dt+d\mu ,\\
&& \dot{p}=-q +c(p-r),\\
&& \dot{r}=-\sigma-c(p-r)+br .
\end{eqnarray*}
From the maximum condition we have
$$
u=\left\{
\begin{array}{cl}
 1, & p+r>0,\\
 -1, & p+r<0.
\end{array}
\right.
$$
Let $y<0$. Then we have the system
\begin{eqnarray*}
&& \dot{q}=a(p-r),\\
&& \dot{\sigma}=-a(p-r),\\
&& \dot{p}=-q +c(p-r),\\
&& \dot{r}=-\sigma-c(p-r)+br ,
\end{eqnarray*}
Since the system is controllable, the subspace $\{ (q,\sigma,p,r)~:~p+r=0\}$ is not invariant. Therefore,  if $y<0$ the optimal control changes its sign at most at a finite number of moments of time.  

\vspace{5mm}

Now suppose that $y(t)=0$, $t\in [t_1,t_2]$, $y(t)<0$, $t<t_1$. If $p+r=0$, then $0=-q-\sigma +br$ and $\dot{q}+\dot{\sigma}=\dot{\mu}$. Hence
$\dot{\mu}=b(-\sigma+2cr+br)$ and we get
$$
\dot{\sigma}=-a(p-r)+b(-\sigma+2cr+br).
$$
Thus for $(q,\sigma,p,r)$ we have the linear system with the matrix
$$
A=\left(
\begin{array}{cccc}
0 & 0 &a &-a\\
0&-b&-a & a+2cb+b^2\\
-1 & 0 &c &-c\\
0&-1&-c&c+b
\end{array}
\right).
$$
The characteristic polynomial of this matrix is 
$$
\lambda^2 (\lambda^2-2c\lambda+2a).
$$
The two-dimensional subspace $L=\{ (q,\sigma,p,r)~:~p+r=0\}$ is now invariant. Obviously $(q,\sigma,p,r)(t)\in L$ only if $p(t_1)+r(t_1)=0$ or $p(t_2)+r(t_2)=0$, i.e., $t_1$ or $t_2$ are the points where the control may change its sign.

\begin{center}
\begin{figure}[h!]
\begin{center}
\pgfimage[width=7cm]{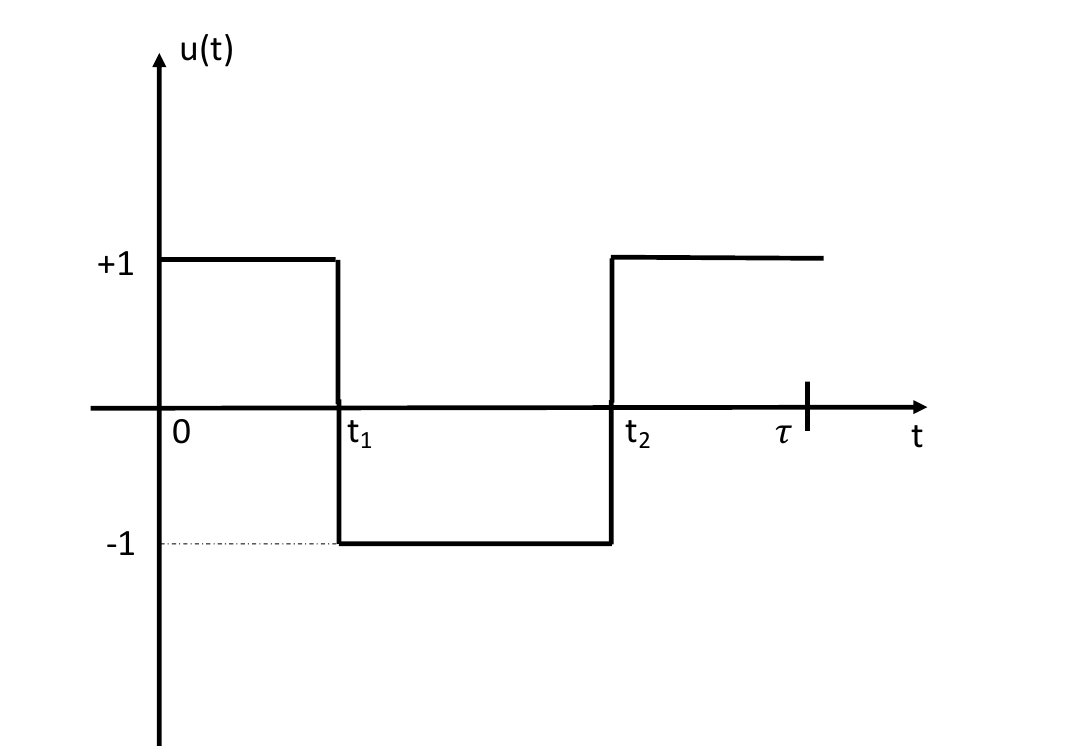}
\end{center}
\caption{\small{ Time-optimal control for Example 2.}}
\label{fig:3}
\end{figure}
\end{center}

The complete analysis of all optimal solutions  in this case is a hard problem. Let us consider a simple, but  interesting, case. 
Assume that  the initial conditions are such that the piston under the control $u=1$ can stop the motion of the cylinder 
at some point $X=Y>0$. Let us imagine for a moment that there is no friction in the system, i.e. $b=c=0$. Then the optimal control has the following structure:
$$
u(t)=
\left\{ 
\begin{array}{rl}
1,& t\in [0,t_1],\\
-1,& t\in [t_1,t_2],\\
1,& t\in [t_2,T],\\
\end{array}
\right.
$$
as illustrated in Figure \ref{fig:3}.

The piston passes the position $Y=0$ at some moment of time $t<t_1$ and  at $t=t_1$ changes the sign {\em before } the collision with the cylinder. At the moment of collision it is already reducing its speed, and after the collision, during the contact with the cylinder, it arrives to the position $Y=0$ with zero speed at $t=t^*$.  Then, it moves to its initial position $Y_0$ as the controlled point body considered in the first example. 

In the presence of negligible  friction, the optimal control driving the system to a small stable neighborhood of  $(y_1,y_1,0,0)$ has the same structure.

\begin{center}
\begin{figure}[h!]
\begin{center}
\pgfimage[width=14cm]{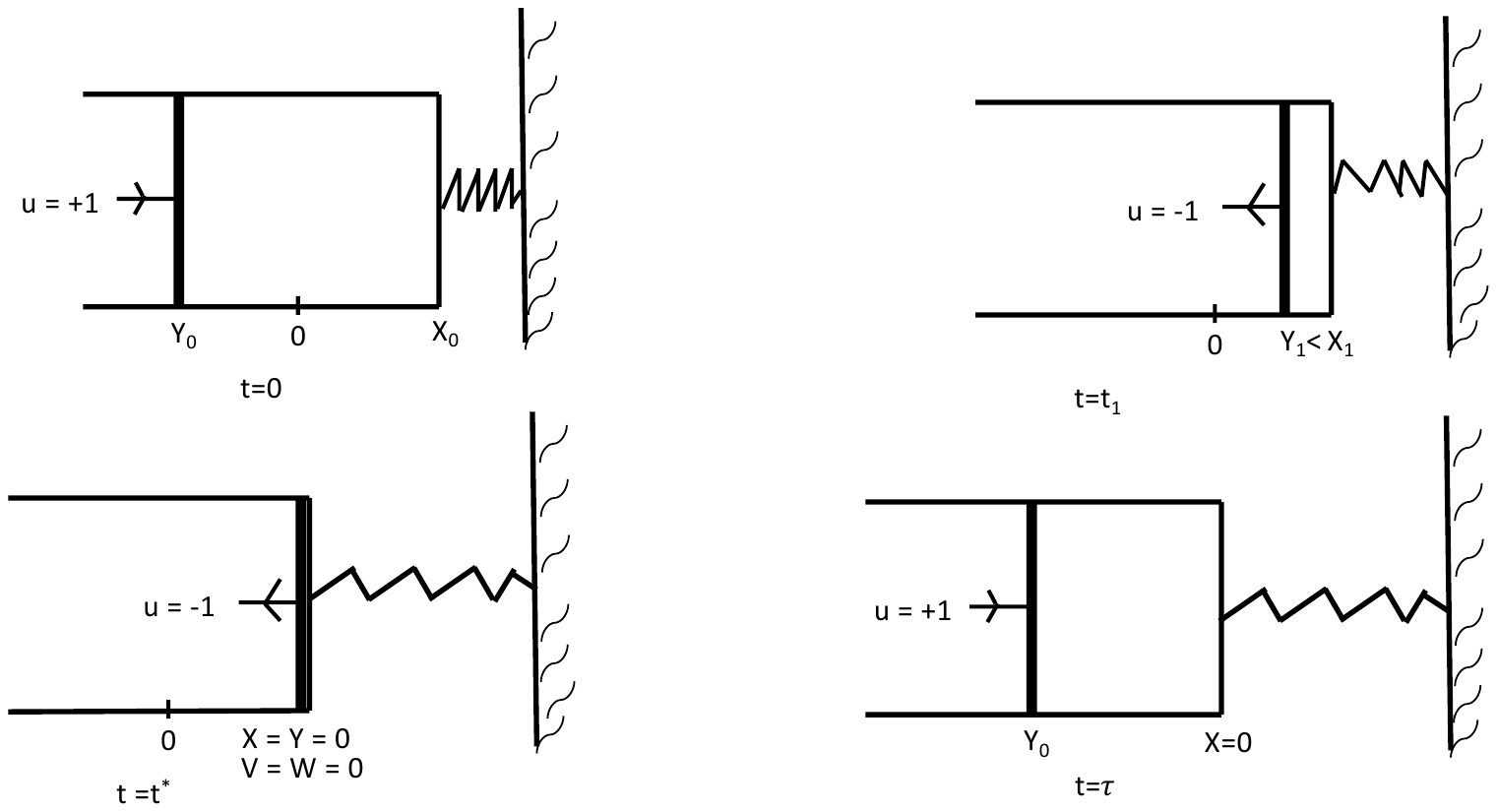}
\end{center}
\caption{\small{The four main phases of the motion (with negligible  friction)  in Example 2. Phase 1, when the cylinder moves away from the wall while  the piston moves in the opposite direction with acceleration, is on the top left figure. Phase 2, on the right top, is when the piston starts braking before collision with the bottom of the cylinder. Phase 3, on the bottom left, is when the piston collides with the bottom of the cylinder. Finally, Phase 4 depicts the situation where the piston returns  to its initial position while the cylinder remains in its equilibrium position.} }
\end{figure}
\end{center}

\newpage

\begin{center}
\textbf{Acknowledgements}
\end{center}

The authors gratefully   thank the support of  Portuguese Foundation for Science and Technology
(FCT) in the framework of the Strategic Funding UIDB/04650\-/2020. 

Also we thank the support  by the ERDF - European Regional Development Fund through the Operational Programme for Competitiveness and Internationalisation - COMPETE 2020, INCO.2030, under the Portugal 2020 Partnership Agreement and by National Funds, Norte 2020, through CCDRN and FCT, within projects \textit{To Chair}  (POCI-01-0145-FEDER-028247), \textit{Upwind} (PTDC/EEI-AUT/31447/2017 - POCI-01-0145-FE\-DER\--03\-14\-47)
and \textit{Systec R\&D unit} (UIDB/00147/2020).

\end{document}